\documentclass[12pt,twoside]{amsart}
\usepackage{amssymb}
\usepackage[all]{xy}

\nonstopmode

\headsep=0.5cm \numberwithin{equation}{section}
\font\tengothic=eufm10 scaled\magstep 1 \font\sevengothic=eufm7
scaled\magstep 1
\newfam\gothicfam
      \textfont\gothicfam=\tengothic
      \scriptfont\gothicfam=\sevengothic

\newtheorem{theorem}{Theorem}[section]

\newtheorem{proposition}[theorem]{Proposition}
\newtheorem{corollary}[theorem]{Corollary}

\theoremstyle{definition}
\newtheorem{definition}[theorem]{Definition} 
\newtheorem{remark}[theorem]{Remark}
\newtheorem{example}[theorem]{Example}

\newtheorem{notation}[theorem]{Notation}
\newtheorem{question}[theorem]{Question}
\newtheorem{problem}[theorem]{Problem}

\newcommand{\Hom}{\operatorname{Hom}}

\newcommand{\Ext}{\operatorname{Ext}}
\newcommand{\rank}{\operatorname{rank}}

\newcommand{\cD}{{\mathcal D}}
\newcommand{\cE}{{\mathcal E}}
\newcommand{\cF}{{\mathcal F}}
\newcommand{\cG}{{\mathcal G}}
\newcommand{\cO}{{\mathcal O}}

\newcommand{\cS}{{\mathcal S}}
\newcommand{\cQ}{{\mathcal Q}}
\newcommand{\cH}{{\mathcal H}}

\newcommand {\FF}{\mathbb{F}}

\newcommand {\CC}{\mathbb{C}}
\newcommand {\ZZ}{\mathbb{Z}}

\newcommand {\PP}{\mathbb{P}}

\def\mapright#1{\smash{ \mathop{\longrightarrow}
    \limits^{#1}}}

\begin{document}
\title[Cohomological characterization of vector bundles] {Cohomological characterization of
vector bundles on multiprojective spaces}

\author[L.\ Costa, R.M.\ Mir\'o-Roig]{L.\ Costa$^*$, R.M.\
Mir\'o-Roig$^{**}$}

\address{Facultat de Matem\`atiques,
Departament d'Algebra i Geometria, Gran Via de les Corts Catalanes
585, 08007 Barcelona, SPAIN } \email{costa@ub.edu}

\address{Facultat de Matem\`atiques,
Departament d'Algebra i Geometria, Gran Via de les Corts Catalanes
585, 08007 Barcelona, SPAIN } \email{miro@ub.edu}

\date{\today}
\thanks{$^*$ Partially supported by BFM2001-3584.}
\thanks{$^{**}$ Partially supported by BFM2001-3584.}



\begin{abstract} We show that Horrock's criterion for the
splitting of vector bundles on $\PP^n$ can be extended to vector
bundles on multiprojective spaces and to smooth projective
varieties with the weak CM property (see Definition
\ref{CMproperty}). As a main tool we use the theory of $n$-blocks
and Beilinson's type spectral sequences. Cohomological
characterizations of vector bundles are also showed.

\end{abstract}


\maketitle

\tableofcontents


 \section{Introduction} \label{intro}

There are two starting points for our work. The first one is the
following  well known result of Horrocks (see \cite{Ho})which
states that a vector bundle on a projective space has no
intermediate cohomology if and only if it decomposes into a direct
sum of line bundles. In \cite{O2}, Ottaviani showed that Horrocks
criterion fails on nonsingular hyperquadrics $Q_3\subset \PP^4$.
Indeed, the Spinor bundle $S$ on $Q_3\subset \PP^4$ has no
intermediate cohomology and it does not decompose into a direct
sum of line bundles. So, it is natural to consider two possible
generalizations of Horrocks criterion  to arbitrary varieties. The
first one consists of characterizing direct sums of line bundles
and the second one consists of characterizing vector bundles
without intermediate cohomology.

Related to the characterization of vector bundles which splits as
direct sum of line bundles; it has been done for vector bundles on
hyperquadrics $Q_n\subset \PP^{n+1}$ and Grassmannians $Gr(k,n)$
by Ottaviani in \cite{O1} and \cite{O2}, respectively. It turns
out that a vector bundle  on $Q_n$ (resp. $Gr(k,n)$) is a direct
sum of line bundles if it has no intermediate cohomology and
satisfies other cohomological conditions involving Spinor bundles
(resp. the tautological $k$-dimensional bundle) and explicitly
written down. Concerning the characterization of vector bundles
without intermediate cohomology besides the result of Horrocks for
vector bundles on projective spaces, there is such a
characterization for vector bundles on hyperquadrics due to
Kn\"{o}rrer; i.e. the line bundles and the Spinor bundles are the
only indecomposable vector bundles on $Q_n\subset \PP^{n+1}$
without intermediate cohomology. Moreover, Buchweitz, Greuel and
Schreyer in \cite{BGS} proved that hyperplanes and hyperquadrics
are the only smooth hypersurfaces in a projective space for which
there are, up to twist, a finite number of indecomposable vector
bundles without intermediate cohomology. See \cite{AG} for the
characterization of vector bundles on $Gr(2,5)$ without
intermediate cohomology and \cite{AC} for the characterization of
rank 2 vector bundles on Fano 3-folds of index 2 without
intermediate cohomology.

The first goal  of this paper is to generalize Horrocks result to
vector bundles on multiprojective spaces $\PP^{n_1}\times \cdots
\times \PP^{n_r}$ and to vector bundles on any smooth projective
variety with the strong CM property (see Definition
\ref{CMproperty}). Indeed, using the notions of exceptional
collections (see Definition \ref{exceptcoll}), $m$-blocks (see
Definition \ref{block}) and the spectral sequences associated to
them (see Theorem \ref{mbe}), we prove that a vector bundle $E $
on $X=\PP^{n_1}\times \cdots \times \PP^{n_r}$  splits provided
$E\otimes \cO _X(t_1,\cdots ,t_r)$ is an ACM bundle for any
$-n_i\le t_i\le 0$, $1\le i \le r$.

\vskip 2mm Our second starting point for this note was another
result of Horrocks which gives a cohomological characterization of
the sheaf of the $p$-differential forms $\Omega ^p_{\PP^n}$ on
$\PP^n$ (\cite{Ho2})and the increasing interest in further
cohomological characterization of vector bundles. Using the notion
of left dual $m$-block collection and again Beilinson's type
spectral sequence, we characterize the $p$-differential forms on
multiprojective spaces.

\vspace{3mm}

Next we outline the structure of this paper. In section 2, we
briefly recall the notions and properties of exceptional sheaf and
full, strongly exceptional collections of sheaves  needed later.
It is well known that the length of any full strongly exceptional
collection of coherent sheaves $\sigma =(E_0,E_1,\cdots ,E_{m})$
on a smooth projective variety $X$ of dimension $n$ is  greater or
equal to $n+1$ and, in \cite{CMR} we call excellent collection any
full exceptional collection of coherent sheaves of length $n+1$.
Excellent collections have nice properties: They are automatically
full strongly exceptional collections and their strong
exceptionality is preserved under mutations. Nevertheless the
existence of an excellent collection on an $n$-dimensional smooth
projective variety  imposes a strong restriction on $X$, namely,
$X$ has to be Fano and $K_0(X)$ a $\ZZ$-free module of rank $n+1$.
In section 3, we generalize the notion of excellent collection
allowing exceptional collections $\sigma =(E_0,E_1,\cdots ,E_{m})$
of arbitrary length but packing the sheaves $E_i$ in suitable
subcollections called blocks.  We introduce the notion of left and
right dual $m$-block collection and we prove its existence
(Proposition \ref{existence-dual}).  In the last part of section
3, we concentrate our attention in varieties $X$ with a number of
blocks generating $D^b(\cO_X-mod)$ one greater than the dimension
of $X$. This leads us to the following definition: We say that an
$n$-dimensional smooth projective variety has the weak CM property
if it has an $n$-block collection which generates
$D^b({\cO}_X$-$mod)$ (see Definition \ref{CMproperty}). Finally,
given a coherent sheaf $\cF$ on a smooth projective variety $X$
with the weak CM property, we derive two Beilinson type spectral
sequences which abuts to $\cF$ (Theorem \ref{mbe}). These two
spectral sequences will play an important role in next section.

 In section 4, we  use
Beilinson type spectral sequence to stablish under which
conditions a vector bundle splits.  As an immediate consequence of
Proposition \ref{main1} we will re-prove: (1) Horrock's criterion
which states that a vector bundle on $\PP^n$ has no intermediate
cohomology if and only if it decomposes into a direct sum of line
bundles (Corollary \ref{coro1}), (2) the characterization  of
vector bundles on a quadric hypersurface $Q_n\subset \PP^{n+1}$,
$n\ge 2$, which splits into a direct sum of line bundles
(Corollary \ref{coro2}) and (3) the characterization of vector
bundles on a Grassmannian $Gr(k,n)$ which splits into a direct sum
of line bundles (Corollary \ref{coro3}). As a main result, we
generalize Horrocks criterion to vector bundles on multiprojective
spaces (see Theorem \ref{mainmulti}) and we get a cohomological
characterization of the $p$-differential forms on multiprojective
spaces (see Theorem \ref{mainmulti2}). We end the paper in \S 5
with some final comments which naturally arise from this paper.
\vspace{3mm}

\vskip 4mm \noindent {\bf Notation} Throughout this paper $X$ will
be a smooth projective variety defined over the complex numbers
$\CC$ and we denote by $\cD=D^b({\cO}_X$-$mod)$  the derived
category of bounded complexes of coherent sheaves of
${\cO}_X$-modules. Notice that $\cD$ is an abelian linear
triangulated category. We identify, as usual, any coherent sheaf
$\cF$ on $X$ to the object $(0 \rightarrow \cF \rightarrow 0) \in
\cD$ concentrated in degree zero and we will not distinguish
between a vector bundle and its locally free sheaf of sections.  A
coherent sheaf $E$ on a smooth projective variety $X$ is an {\bf
ACM sheaf} if $H^{i}(X,E\otimes \cO_X(t))=0$ for any $i$,
$0<i<\dim X$, and for any $t\in \ZZ$; and we say that $E$ has {\bf
no intermediate cohomology} if and only if $H^{i}(X,E\otimes L)=0$
for any $i$, $0<i<\dim X$, and for any line bundle $L$ on $X$.


\section{Preliminaries}

As we pointed out in the introduction, in this section we gather
the basic definitions and properties on exceptional sheaves,
exceptional collections of sheaves, strongly exceptional
collections of sheaves and full exceptional collections of sheaves
needed in the sequel.

\begin{definition}\label{exceptcoll}
Let $X$ be a smooth projective variety.

(i) An object  $F \in \cD$ is {\bf exceptional} if
$\Hom^{\bullet}_{\cD} (F,F)$ is a 1-dimensional algebra generated
by the identity.

(ii) An ordered collection $(F_0,F_1,\cdots ,F_m)$  of objects of
$\cD$ is an {\bf exceptional collection} if each object $F_{i}$ is
exceptional and $\Ext^{\bullet}_{\cD}(F_{k},F_{j})=0$ for $j<k$.

(iii) An exceptional collection $(F_0,F_1,\cdots ,F_m)$ of objects
of $\cD$ is a {\bf strongly exceptional collection} if in addition
$\Ext^{i}_{\cD}(F_j,F_k)=0$ for $i\neq 0$ and  $j \leq k$.

(iv) An ordered collection of objects of $\cD$, $(F_0,F_1,\cdots
,F_m)$,  is a {\bf full (strongly) exceptional collection} if it
is a (strongly) exceptional collection  and $F_0$, $F_1$, $\cdots
$ , $F_m$ generate the bounded derived category $\cD$.
\end{definition}

\begin{remark} The existence of a full strongly
exceptional collection $(F_0,F_1,\cdots,F_m)$ of coherent sheaves
on a smooth projective variety $X$ imposes rather a strong
restriction on $X$, namely that the Grothendieck group
$K_0(X)=K_0(\cO _X-mod)$ is isomorphic to $\ZZ^{m+1}$.
\end{remark}

\begin{example} \label{prihirse2} (1) ($\cO_{\PP^r}(-r)$, $\cO_{\PP^r} (-r+1) $,
 $\cO_{\PP^r} (-r+2)
$, $\cdots $, $\cO_{\PP^r}$) is a full strongly exceptional
collection of coherent sheaves on a projective space $\PP^r$ and
($\cO_{\PP^r}$, $\Omega^1_{\PP^r} (1) $, $\Omega^2_{\PP^r} (2) $,
$\cdots $, $\Omega^r_{\PP^r} (r) $) is also a full strongly
exceptional collection of coherent sheaves on $\PP^r$.

(2) Let $\FF_n=\PP(\cO_{\PP^1}\oplus \cO_{\PP^1}(n))$, $n\ge 0$,
be a Hirzebruch surface. Denote by $\xi$ (resp. $F$) the class of
the tautological line bundle (resp. the class of a fiber of the
natural projection $p: \FF_n \rightarrow \PP^1$). Then, ($\cO$,
$\cO (F) $, $\cO (\xi) $, $\cO (F+\xi) $) is a full strongly
exceptional collection of coherent sheaves on $\FF_n$.

(3) Let $\pi: \widetilde{\PP}^2(l) \rightarrow \PP^2$  be the blow
up of $\PP^2$ at $l$ points and let $L_1=\pi^{-1}(p_1),
\ldots,L_l=\pi^{-1}(p_l)$ be the exceptional divisors. Then,
\[ (\cO, \cO(L_1), \cO(L_2), \cdots, \cO(L_l), \cO(H),\cO(2H)) \]
 is a full strongly exceptional collection of coherent sheaves on $\widetilde{\PP}^2(l)$.

(4) Let $\cE$ be a rank $r$ vector bundle on a smooth projective
variety $X$. If $X$ has a full strongly exceptional collection of
line bundles then $\PP(\cE)$ also has a full strongly exceptional
collection of line bundles. In particular, any $d$-dimensional,
smooth, complete toric variety $V$ with a splitting fan $\Sigma
(V)$ has a full strongly exceptional collection of line bundles
and any $d$-dimensional, smooth, complete toric variety $V$ with
Picard number 2 or, equivalently, with $d+2$ generators has a full
strongly exceptional collection of line bundles (see \cite{CMR}).

(5) $( \cO_{\PP^n}(-n) \boxtimes \cO_{\PP^m}(-m),
\cO_{\PP^n}(-n+1) \boxtimes \cO_{\PP^m}(-m),
 \cdots, \cO_{\PP^n} \boxtimes \cO_{\PP^m}(-m),
\cdots,  \cO_{\PP^n}(-n) \boxtimes \cO_{\PP^m}, \cO_{\PP^n}(-n+1)
\boxtimes \cO_{\PP^m},$ $ \cdots, \cO_{\PP^n} \boxtimes
\cO_{\PP^m}) $ is a full strongly exceptional collection of
locally free sheaves on $\PP^n \times \PP^m$.
\end{example}

\vspace{3mm}

We have seen many examples of smooth projective varieties which
have a full strongly exceptional collection of line bundles and we
want to point out that there are many other examples of smooth
projective varieties which have a full strongly exceptional
collection of bundles of higher rank but they  don't have a full
strongly exceptional collection of line bundles.

 \vspace{3mm}

\begin{example}
\label{prihirse}
 (1) Let  $X=Gr(k,n)$ be the Grassmannian of
$k$-dimensional subspaces of the $n$-dimensional vector space.
Assume $k>1$. We have $Pic(X) \cong \ZZ \cong < \cO_{X}(1)>$, $K_X
\cong \cO_X(-n)$ and the canonical exact sequence
\[ 0 \rightarrow \cS \rightarrow \cO_X^n \rightarrow \cQ \rightarrow 0
\] where $\cS$ denotes the tautological $k$-dimensional bundle and
$\cQ$ the quotient bundle.

In the sequel, $\Sigma^{\alpha}\cS$ denotes the space of the
irreducible representations of the group $GL(\cS)$ with highest
weight $\alpha=(\alpha_1, \ldots, \alpha_s)$ and
$|\alpha|=\sum_{i=1}^{s} \alpha_i$. Denote by $A(k,n)$ the set of
locally free sheaves $\Sigma^{\alpha}\cS$ on $Gr(k,n)$ where
$\alpha$ runs over Young diagrams fitting inside a $k \times
(n-k)$ rectangle. Set $\rho(k,n):= \sharp A(k,n)$. By \cite{Kap};
Proposition 2.2 (a) and Proposition 1.4, $A(k,n)$ can be totally
ordered in such a way that we obtain a full strongly exceptional
collection ($ E_1, \ldots, E_{\rho(k,n)}$) of locally free sheaves
on $X$. Notice that $\cS \in A(k,n)$ has rank $k$ and hence this
collection has locally free sheaves of rank greater than one. In
addition, any full strongly exceptional collection of coherent
sheaves on $X$ has a sheaf of rank greater than one. Indeed,  any
full strongly exceptional collection of coherent sheaves on $X$
has the same length  equals to the rank $\rho(k,n)$ of the
Grothendieck group of $X$. On the other hand, since $Pic(X) \cong
< \cO_{X}(1)>$ and $K_X \cong \cO_X(-n)$, any full strongly
exceptional collection of coherent sheaves has at most $n+1$
summands which are line bundles. Therefore, since $n+1 <
\rho(k,n)=rk(K_0(X))$, any full strongly exceptional collection
has a sheaf of rank different from one.

(2)  Any full strongly exceptional collection of locally free
sheaves on a hyperquadric  $Q_n \subset \PP^{n+1}$, $n>2$, has a
sheaf of rank different from one. In fact, if $n \geq 3$ then
Pic$(Q_n)=\ZZ =\langle \cO_{Q_n}(1) \rangle$, $K_{Q_n} \cong
\cO_{Q_n}(-n)$ and
\[ rank(K_0(Q_n))= \begin{cases} n+1 \quad \mbox{if } n \quad \mbox{is
odd} \\
 n+2 \quad \mbox{if } n \quad \mbox{is
even.}  \end{cases} \] Moreover, by \cite{Kap2}; Proposition 4.9,
if $n$ is even and $\Sigma _1$, $\Sigma _2$ are the Spinor bundles
on $Q_n$, then
\[(\Sigma _1(-n), \Sigma _2(-n), \cO_{Q_n}(-n+1), \cdots, \cO_{Q_n}(-1), \cO_{Q_n}) \]
is a full strongly exceptional collection of locally free sheaves
on $Q_n$; and if $n$ is odd and $\Sigma$ is the Spinor bundle on
$Q_n$, then
\[( \Sigma(-n), \cO_{Q_n}(-n+1), \cdots, \cO_{Q_n}(-1), \cO_{Q_n}) \]
is a full strongly exceptional collection of locally free sheaves
on $Q_n$.
\end{example}

 \vspace{3mm}

\begin{definition}
\label{mutation} Let $X$ be a smooth projective variety and let
$(A,B)$ be an exceptional pair of objects of $\cD$. We define
objects $L_AB$ and $R_BA$ with the aid of the following
distinguished triangles in the category $\cD$:
\begin{equation}  \label{t1} L_AB \rightarrow \Hom^{\bullet}_{\cD}(A,B) \otimes A \rightarrow B
\rightarrow L_AB[1]
\end{equation}
\begin{equation}  \label{t2} R_BA[-1] \rightarrow A \rightarrow \Hom^{\times \bullet}_{\cD}(A,B)
 \otimes B \rightarrow R_BA.
\end{equation}
\end{definition}

\vspace{2mm}

\begin{notation} \label{composicio} Let $X$ be a smooth projective variety and let $\sigma=(F_0,
\cdots, F_m )$ be  an exceptional collection of objects of $\cD$.
It is convenient to agree that for any $0 \leq i,j \leq m$ and
$i+j \leq m$,
\[  R^{(j)}F_{i}=  R^{(j-1)}RF_{i}= R_{F_{i+j}} \cdots \cdots
  R_{F_{i+2}} R_{F_{i+1}} F_{i}=: R_{F_{i+j} \cdots \cdots
  F_{i+2} F_{i+1}} F_{i} \]
 and similar notation for compositions of left mutations.
\end{notation}

\vspace{3mm}

If $X$ is a smooth projective variety and  $\sigma=(F_0, \cdots,
F_m )$ is  an exceptional collection of objects of $\cD$, then any
mutation of $\sigma$ is an exceptional collection. Moreover, if
$\sigma$ generates the category $\cD$, then the mutated collection
also generates $\cD$.

 Nevertheless, in general, a mutation of a strongly exceptional collection is
not a  strongly exceptional collection. In fact, take $X= \PP^1
\times \PP^1$ and consider the full strongly exceptional
collection $\sigma=(\cO_X, \cO_X(1,0), \cO_X(0,1), \cO_X(1,1) ) $
of line bundles on $X$. It is not difficult to check that the
mutated collection $$ (\cO_X, \cO_X(1,0),
L_{\cO_X(0,1)}\cO_X(1,1), \cO_X(0,1) )=(\cO_X, \cO_X(1,0),
\cO_X(-1,1), \cO_X(0,1) )$$ is no more a strongly exceptional
collection of line bundles on $X$.


\section{$m$-blocks and Beilinson's spectral sequence}

Let $X$ be a smooth projective variety of dimension $n$. It is
well known that all full strongly exceptional collections of
coherent sheaves on $X$ have the same length and it is equal to
the rank of $K_0(X)$. Even more, this length is bounded below by
$n+1$ because for any smooth projective variety $X$ of dimension
$n$  we have $rank(K_0(X)) \geq n+1$. In \cite{CMR3}; we give the
following definition (see also \cite{BP} and \cite{H}).

\vspace{3mm}

\begin{definition}
\label{defexcellent}  Let $X$ be a smooth projective variety of
dimension $n$. We say that an ordered collection of coherent
sheaves $\sigma=(E_0, \cdots, E_n)$ is an {\bf excellent
collection} if it is a full exceptional collection of coherent
sheaves on $X$ of minimal length, $n+1$, i.e. of length one
greater than the dimension of $X$.
\end{definition}

\vspace{3mm}

By  \cite{Bo}; Assertion 9.2, Theorem 9.3 and Corollary 9.4,
excellent collections are automatically strongly exceptional
collections of coherent sheaves and  the strongly exceptionality
is preserved under mutations.

\begin{example}
\label{exemplesexcellents} (1)  The collection
$\sigma=(\cO_{\PP^r}(-r), \cO_{\PP^r} (-r+1) , \cO_{\PP^r} (-r+2)
, \cdots , \cO_{\PP^r} )$ of line bundles on $\PP^r$ is an
excellent collection of coherent sheaves.

(2) If $n$ is odd and $Q_n \subset \PP^{n+1}$ is a quadric
hypersurface,  the collection of locally free sheaves
\[( \Sigma(-n), \cO_{Q_n}(-n+1), \cdots, \cO_{Q_n}(-1), \cO_{Q_n}) \]
being $\Sigma$ the Spinor bundle on $Q_n$ is an excellent
collection of locally free sheaves on $Q_n$.

(3) If $n$ is even and $Q_n \subset \PP^{n+1}$ is a quadric
hypersurface,  the collection of locally free sheaves
\[(\Sigma _1(-n), \Sigma _2(-n), \cO_{Q_n}(-n+1), \cdots, \cO_{Q_n}(-1), \cO_{Q_n}) \]
being $\Sigma _1$ and $\Sigma _2$ the Spinor bundles on $Q_n$, is
a full strongly exceptional collection of locally free sheaves on
$Q_n$. Since all full strongly exceptional collections of coherent
sheaves on $Q_n$ have  length  $n+2$ we conclude that there are no
excellent collections of coherent sheaves on $Q_n$ for even $n$.

(4) It follows from Example \ref{prihirse} that there are no
excellent collections of coherent sheaves on $Gr(k,n)$ if $k \neq
n-1$.

(5) Any smooth Fano threefold $X$ with $Pic(X) \cong \ZZ$ and
trivial intermediate Jacobian has an excellent collection (see
\cite{CMR3}; Proposition 3.6).
\end{example}

\vskip 2mm It is an interesting problem to characterize the smooth
projective varieties which have an excellent collection. We want
to stress that the existence of an excellent collection on an
$n$-dimensional smooth variety $X$ imposes a strong restriction on
$X$; e.g. $X$ has to be a Fano variety (\cite{BP}; Theorem 3.4)
and the Grothendieck group $K_0(X)$ has to be a $\ZZ $-free module
of rank $n+1$. So, it is convenient to generalize the notion of
excellent collection in order to be able to apply the results
derived from its existence to varieties as Grassmannians,
even-dimensional hyperquadrics, multiprojective spaces, etc, which
do not have excellent collections. This will be achieved allowing
exceptional collections $\sigma=(F_0, \ldots,F_m)$ of arbitrary
length but packing the objects $F_i \in \cD$ in suitable
subcollections called blocks. The notion of block was introduced
by Karpov and Nogin in \cite{KN} and we will recall its definition
and properties (see also \cite{H}).

\vspace{3mm}

\begin{definition}
\label{block} (i) An exceptional collection $(F_0,F_1, \cdots,
F_m)$ of objects of $\cD$ is a {\bf block} if
$\Ext^{i}_{\cD}(F_j,F_k)=0$ for any $i$ and  $j \neq k$.

(ii) An {\bf $m$-block collection of  type} $(\alpha_0, \alpha_1,
\cdots, \alpha_m)$ of objects of $\cD$ is an exceptional
collection
\[(\cE_0, \cE_1, \cdots, \cE_m)=(E_1^0,\cdots,E_{\alpha_0}^0, E_1^1,
 \cdots, E_{\alpha_1}^1,
\cdots, E_1^m,\cdots, E_{\alpha_m}^m) \] such that all the
subcollections $\cE_i=(E_1^i,E_2^i,\cdots, E_{\alpha_i}^i)$ are
blocks.
\end{definition}

Note that an exceptional collection $(E_0,E_1, \cdots, E_m)$ is an
$m$-block of type $(1,1, \cdots, 1)$.

\vspace{3mm}

\begin{example}
\label{exempleblock} (1)  Let  $X=Gr(k,n)$ be the Grassmannian of
$k$-dimensional subspaces of the $n$-dimensional vector space,
$k>1$. In Example \ref{prihirse} (1), we have seen that $A(k,n)$
can be totally ordered in such a way that we obtain a full
strongly exceptional collection \[ \sigma=( E_1, \ldots,
E_{\rho(k,n)})\] of locally free sheaves on $X$. On the other
hand, by \cite{Kap2} (3.5), $\Hom(\Sigma^{\alpha} \cS,
\Sigma^{\beta}\cS) \neq 0$ only if $\alpha_i \geq \beta_i$ for all
$i$. So, packing in the same block $\cE _r$ the bundles
$\Sigma^{\alpha} \cS \in \sigma$ with $|\alpha|= k(n-k)-r$ and
taking into account that $0 \leq |\alpha| \leq k(n-k)$ we obtain
\[ \sigma=( E_1, \ldots, E_{\rho(k,n)})=(\cE_0, \ldots, \cE_{k(n-k)}) \]
a $k(n-k)$-block collection of vector bundles on $X$.

\vskip 2mm (2) Let $Q_n \subset \PP^{n+1}$, $n \geq 2$, be a
hyperquadric variety. According to Example \ref{prihirse} (2), if
$n$ is even and $\Sigma _1$, $\Sigma _2$ are the Spinor bundles on
$Q_n$, then
\[(\Sigma _1(-n), \Sigma _2(-n), \cO_{Q_n}(-n+1), \cdots, \cO_{Q_n}(-1), \cO_{Q_n}) \]
is a full strongly exceptional collection of locally free sheaves
on $Q_n$; and if $n$ is odd and $\Sigma$ is the Spinor bundle on
$Q_n$, then
\[( \Sigma(-n), \cO_{Q_n}(-n+1), \cdots, \cO_{Q_n}(-1), \cO_{Q_n}) \]
is a full strongly exceptional collection of locally free sheaves
on $Q_n$. Since $\Ext^i(\Sigma_1,\Sigma_2)=0$ for any $i \geq 0$,
we get that $(\cE_0, \cE_1, \ldots, \cE_n)$ where  \[
\cE_i=\cO_{Q_n}(-n+i) \quad \mbox{for } 1 \leq i \leq n, \quad
\cE_0=
\begin{cases} (\Sigma_1(-n), \Sigma_2(-n))  & \mbox{if } n
\quad \mbox{even}
\\ (\Sigma(-n)) & \mbox{if } n \quad \mbox{odd} \\
\end{cases} \] is an $n$-block collection of coherent sheaves on $Q_n$ for all
$n$. \vskip 2mm

(3) Let $X=\PP^{n_1} \times \cdots \times \PP^{n_s}$ be a
multiprojective space of dimension $d=n_1+\cdots+n_s$. For any $1
\leq i \leq s$, denote by $p_i:X \rightarrow \PP^{n_i}$ the
natural projection and write
 \[ \cO_X(a_1, a_2, \cdots, a_s):= p_1^* \cO_{\PP^{n_1}}(a_1)
 \otimes p_2^*\cO_{\PP^{n_2}}(a_2) \otimes \cdots \otimes p_s^* \cO_{\PP^{n_s}}(a_s).\]
 For any $0 \leq j \leq d$, denote by $\cE_j$
 the collection of all line bundles on $X$
 \[ \cO_X(a_1^j, a_2^j, \cdots, a_s^j)\]
 with $-n_i \leq a_i^j \leq 0$ and
 $\sum_{i=1}^{s}a_i^j=j-d$.  Using the K\"{u}nneth formula for
 locally free sheaves on algebraic varieties, we prove that
 each $\cE_j$ is a block and that
 \[(\cE_0, \cE_1, \cdots, \cE_{d}) \]
 is a $d$-block collection of line bundles
 on $X$.
\end{example}

\vspace{3mm}

We will now introduce the notion of mutation of block collections.

\vspace{3mm}

\begin{definition} \label{mutationblock}
Let $X$ be a smooth projective variety and consider
  a 1-block
collection $(\cE,\cF)=(E_1, \cdots, E_n,F_1, \cdots, F_m)$ of
objects of $\cD$. A {\bf left mutation} of $F_j$ by $\cE$ is the
object defined by (see Notation \ref{composicio})
\[L_{\cE}F_j:=L_{E_1E_2 \cdots E_n}F_j \]
and a  {\bf right mutation} of $E_j$ by $\cF$ is the object
defined by
\[R_{\cF}E_j:= R_{F_mF_{m-1} \cdots F_1}E_j. \]
A {\bf left mutation} of $(\cE, \cF)$ is the pair $(L_{\cE}\cF,
\cE)$ where
\[L_{\cE}\cF:=(L_{\cE}F_1, L_{\cE}F_2, \cdots,L_{\cE}F_m) \]
and a {\bf right mutation} of $(\cE, \cF)$ is the pair
$(\cF,R_{\cF}\cE)$ where
\[R_{\cF}\cE:=(R_{\cF}E_1, R_{\cF}E_2, \cdots, R_{\cF}E_n). \]
\end{definition}

 \vspace{3mm}

\begin{remark} By \cite{GK} (2.2), for any exceptional object $X
\in\cD$, any pair of object $F, G \in \cD$ and any integer $i$ we
have:
\[ Ext^i_{\cD}(L_XF,L_XG)= Ext^i_{\cD}(F,G),\]
\[Ext^i_{\cD}(R_XF,R_XG)= Ext^i_{\cD}(F,G). \]
Hence, for any 1-block collection $(\cE, \cF)=(E_1, \ldots,
E_n,F_1 \ldots, F_m )$ and integers $j \neq k$,
\[ Ext^i_{\cD}(L_{\cE}F_j,L_{\cE}F_k)=Ext^i_{\cD}(L_{E_1 \cdots E_n}F_j,L_{E_1 \cdots E_n}F_k)
=Ext^i_{\cD}(F_j,F_k),\]
\[ Ext^i_{\cD}(R_{\cF}E_j,R_{\cF}E_k)= Ext^i_{\cD}(R_{F_m \cdots F_1}E_j,R_{F_m \cdots F_1}E_k)
=Ext^i_{\cD}(E_j,E_k) \] and thus both $L_{\cE}\cF$ and
$R_{\cF}\cE$ are blocks and the pairs $(L_{\cE}\cF,\cE)$ and
$(\cF, R_{\cF}\cE)$ are 1-block collections.
\end{remark}

\begin{remark} \label{triangles}
It follows from the proof of \cite{KN}; Proposition 2.2 and
Proposition 2.3 that given a 1-block collection $(\cE,\cF)=(E_1,
\cdots, E_n,F_1, \cdots, F_m)$, the objects $L_{\cE}F_j$ and
$R_{\cF}E_j$ can be defined with the aid of the following
distinguished triangles in the category $\cD$
\begin{equation}
\label{triangleleft} L_{\cE}F_j \rightarrow \oplus_{i=1}^{n}
\Hom^{\bullet}_{\cD}(E_i,F_j) \otimes E_i \rightarrow F_j
\rightarrow L_{\cE}F_j[1]
\end{equation}
\begin{equation}
\label{triangleright} R_{\cF}E_j[-1] \rightarrow E_j \rightarrow
\oplus_{i=1}^{m} \Hom^{\times \bullet}_{\cD}(E_j,F_i)\otimes F_i
\rightarrow R_{\cF}E_j.
\end{equation}
\end{remark}

\vspace{3mm}

Applying $\Hom^{\bullet}_{\cD}(E_i, *)$ to the triangle
(\ref{triangleleft}) we get the orthogonality relation

\begin{equation}
\label{orto1l} \Hom^{\bullet}_{\cD}(E_i, L_{\cE} F_j)=0 \quad
\mbox{for all} \quad 1 \leq i \leq n
\end{equation}

\vskip 2mm \noindent i.e., $L_{\cE}F_j \in [\cE]^{\bot}:= \{F \in
\cD | \Hom^{\bullet}_{\cD}(E,F)=0 \quad \mbox{for all } E \in
[\cE]\}$, where we denote by $[\cE]$ the full triangulated
subcategory of $\cD$ generated by $E_1, \cdots, E_n$.

Similarly, $\Hom^{\bullet}_{\cD}(*, F_j)$ applied to the triangle
(\ref{triangleright}) gives the orthogonality relation

\begin{equation}
\label{orto1r} \Hom^{\bullet}_{\cD}( R_{\cF} E_i, F_j)=0 \quad
\mbox{for all} \quad 1 \leq j \leq m
\end{equation}

\vskip 2mm \noindent i.e., $R_{\cF}E_i \in$  $^{\bot}[\cF]:= \{E
\in \cD | \Hom^{\bullet}_{\cD}(E,F)=0 \quad \mbox{for all } F \in
[\cF]\}$.

Taking $E' \in$ $^{\bot}[\cF]$ and $E'' \in [\cE]^{\bot}$ and
applying $\Hom^{\bullet}_{\cD}(E', *)$ and
$\Hom^{\bullet}_{\cD}(*, E'')$ to the triangles
(\ref{triangleleft}) and (\ref{triangleright}) we get for any $H
\in \cD$
\begin{equation}
\label{orto2l} \Hom^{\bullet}_{\cD}(E',H)=\Hom^{\bullet}_{\cD}(E',
L_{\cE}H)[1],
\end{equation}
\begin{equation}
 \label{orto2r} \Hom^{\bullet}_{\cD}(H,
 E'')=\Hom^{\bullet}_{\cD}(R_{\cE}H,E'')[1].
\end{equation}

  \vspace{3mm}

\begin{notation} It is convenient to agree that
 \[  R^{(j)}\cE_{i}=  R^{(j-1)}R\cE_{i}= R_{\cE_{i+j}} \cdots \cdots
  R_{\cE_{i+2}} R_{\cE_{i+1}} \cE_{i}=: R_{\cE_{i+j} \cdots \cdots
  \cE_{i+2} \cE_{i+1}} \cE_{i} \]
   \[  L^{(j)}\cE_{i}=  L^{(j-1)}L\cE_{i}= L_{\cE_{i-j}} \cdots \cdots
  L_{\cE_{i-2}} L_{\cE_{i-1}} \cE_{i}=: L_{\cE_{i-j} \cdots \cdots
  \cE_{i-2} \cE_{i-1}} \cE_{i}. \]
\end{notation}

\vspace{3mm}

Let $\sigma=(\cE_0, \cdots, \cE_m)$ be an $m$-block collection of
type $\alpha_0, \cdots, \alpha_m$ of objects of $\cD$ which
generates $\cD$. Two $m$-block collections $\cH=(\cH_0, \cdots,
\cH_m)$ and $\cG=(\cG_0, \cdots, \cG_m)$ of type  $\beta_0,
\cdots, \beta_m$ with $\beta_i=\alpha_{m-i}$ of objects of $\cD$
are called {\bf left dual $m$-block collection of $\sigma $} and
{\bf right dual $m$-block collection of $\sigma $} if
\begin{equation}\label{orto1}  \Hom^{\bullet}_{\cD}(H_{j}^{i},E^k_l)=\Hom^{\bullet}_{\cD}(E^k_l,
G^i_j)=0 \end{equation} except for

\begin{equation}\label{orto2}
 \Ext^{k}_{\cD}(H_{i}^{k},E_i^{m-k})=\Ext^{m-k}_{\cD}(E_i^{m-k},
G_i^k)=\CC.
\end{equation}
\vspace{3mm}

\begin{proposition}\label{existence-dual}
\label{dual} Left dual $m$-block collections and right dual
$m$-block collections exist and they are unique up to isomorphism.
\end{proposition}
\begin{proof} Let $\sigma =(\cE_0, \cdots, \cE_m)$ be an $m$-block
collection of type type $\alpha_0, \cdots, \alpha_m$ of objects of
$\cD$. We will construct explicitly the left and the right dual
$m$-block collection of $\sigma $ by consequent mutations of the
$m$-block collection $\sigma $. We consider \begin{equation}
\label{caracdelsr} \cH=(R^{(0)}\cE_m, R^{(1)}\cE_{m-1}, \cdots,
R^{(m)}\cE_0) \end{equation} where by definition
\[ \begin{array}{rl} R^{(i)}\cE_{m-i}= & (R^{(i)}E_1^{m-i}, \cdots,
R^{(i)}E_{\alpha_{m-i}}^{m-i} ) \\ = & (R_{\cE_m\cE_{m-1} \cdots
\cE_{m-i+1} }E_1^{m-i}, \cdots, R_{\cE_m\cE_{m-1} \cdots
\cE_{m-i+1} }E_{\alpha_{m-i}}^{m-i} ). \end{array} \]

Let us check that it satisfies the orthogonality conditions
$(\ref{orto1})$ and $(\ref{orto2})$. It follows from
(\ref{orto1r}) that $R_{\cE_m\cE_{m-1} \cdots \cE_{m-i+1}
}E_{k}^{m-i} \in$ $^{\bot}[\cE_{m-i+1}, \cdots, \cE_{m}]$ and
hence for any $l$ with $m-i+1 \leq l \leq m $ and any $j$ with $1
\leq j \leq \alpha_l$
\[  \Hom^{\bullet}_{\cD}(R_{\cE_m\cE_{m-1} \cdots
\cE_{m-i+1} }E_{k}^{m-i}, E_j^l)=0. \]

On the other hand, since $\sigma $ is an exceptional collection,
for any $l$ with $0 \leq l \leq m-i$, and any $p$ with $m-i+1 \leq
p \leq m$
\[ \Hom^{\bullet}_{\cD}(E_q^p, E_j^l)=0, \quad 1 \leq q \leq \alpha_p, \quad 1 \leq j \leq \alpha_l. \]
So, for any $l$ with $0 \leq l \leq m-i$ and any $j$ with $1 \leq
j \leq \alpha_l$, $E_j^l \in$ $^{\bot}[\cE_{m-i+1}, \cdots,
\cE_m]$ and applying repeatedly (\ref{orto2r}) we get
\[  \Hom^{\bullet}_{\cD}(R_{\cE_m\cE_{m-1} \cdots
\cE_{m-i+1} }E_{k}^{m-i}, E_j^l)= \Hom^{\bullet}_{\cD}(E_k^{m-i},
E_j^l)[-i] = \begin{cases}  0  \mbox{ if } l < m-i \\ \CC \mbox{
in degree $i$ if } l=m-i. \end{cases} \] Therefore, $\cH$ is
indeed the left dual $m$-block collection of $\sigma $. By
consequent left mutations of the $m$-block collection $\sigma $
and arguing in the same way we get the right dual $m$-block
collection of $\sigma $.
\end{proof}

\vspace{3mm}

We want to point out that the notion of $m$-block collection is
the convenient generalization of the notion of excellent
collection we were looking for. Indeed, we will see that the
behavior of $n$-block collections, $n=\dim(X)$, is really good in
the sense that they are automatically strongly exceptional
collections and that their structure is preserved under mutations
through blocks. More precisely we have:

\begin{proposition}
\label{bonespropietats} Let $X$ be a smooth projective variety of
dimension $n$ and let $\sigma=(\cE_0, \cdots, \cE_n)$ be an
$n$-block collection of coherent sheaves on $X$ and assume that
$\sigma$ generates de category $\cD$. Then we get:

(1) The sequence $\sigma$ is a full strongly exceptional
collection of coherent sheaves on $X$.

(2) All mutations through the blocks $\cE_i$ can be computed using
short exact sequences of coherent sheaves.

(3) Any mutation of $\sigma$ through any block $\cE_i$ is a full
strongly exceptional collection of pure sheaves, i.e. complexes
concentrated in the zero component of the grading.

(4) Any mutation of $\sigma$ through any block $\cE_i$ is an
$n$-block collection.
\end{proposition}
\begin{proof} See \cite{Bo}; Theorem 9.5 and Remark b) below and
\cite{H}; Theorem 1.
\end{proof}

\vspace{3mm}

These nice properties  led us to introduce the following
definition

\vspace{3mm}

\begin{definition}
\label{CMproperty} Let $X$ be a smooth projective variety of
dimension $n$. We say that $X$ has the {\bf weak CM property} if
there exists an $n$-block collection $(\cE_0, \cdots, \cE_n)$ of
type $(\alpha_0, \ldots, \alpha_n)$ of coherent sheaves  on $X$
which generates $\cD$. We say that $X$ has the {\bf  CM property}
if in addition, for all $E_{i}^{n}\in \cE_{n}$ and all $E_{l}^k\in
\cE_k$ with $0 \leq k \leq n-1$, $E_{i}^n\otimes E_{l}^k$ is an
ACM sheaf; and finally we say that $X$ has the {\bf strong CM
property} if in addition, all the exceptional coherent sheaves
$E_j^i \in \cE_i$ are line bundles.
\end{definition}

\vspace{3mm}

\begin{remark}
We want to point out that the number of blocks is one greater than
the dimension of $X$ but a priori there is no restriction on the
length $\alpha_i$ of each block $\cE_i=(E_1^i,
\ldots,E_{\alpha_i}^i)$.
\end{remark}

\vspace{3mm}

It is clear that any smooth projective variety with an excellent
collection has the weak CM property. Let us now see many examples
of varieties with the (weak) CM property which do not have
excellent collections of coherent sheaves.

\begin{example}
\label{CMexample}

(1) Since any line bundle on $\PP^n$ is ACM, it follows from
Example \ref{exemplesexcellents} (1) that $\PP^n$ has the strong
CM property.

(2) Let $Q_n \subset \PP^{n+1}$, $n\ge 2$, be a hyperquadric
variety. According to Example \ref{exempleblock} (2), $
\sigma=(\cE_0, \cE_1, \ldots, \cE_n)$ where  \[
\cE_i=\cO_{Q_n}(-n+i) \quad \mbox{for } 1 \leq i \leq n, \quad
\cE_0=
\begin{cases} (\Sigma_1(-n), \Sigma_2(-n))  & \mbox{if } n
\quad \mbox{even}
\\ (\Sigma(-n)) & \mbox{if } n \quad \mbox{odd} \\
\end{cases} \] is an $n$-block collection of coherent sheaves on $Q_n$ for all
$n$. Since spinor bundles and line bundles on $Q_n$ are ACM
bundles and $\cE_n=\cO_{Q_n}$, we deduce that $Q_n$ has the CM
property.

(3) Let $X=\PP^{n_1} \times \cdots \times \PP^{n_s}$ be any
multiprojective space and let $\sigma=(\cE_0, \cdots, \cE_{n_1 +
\cdots + n_s })$ be the $(n_1+\cdots+n_s)$-block collection of
line bundles on $X$ given in Example \ref{exempleblock} (3). Using
the K\"{u}nneth formula, the fact that $H^{\alpha}(\PP^{n_j},
\cO_{\PP^{n_j}}(a))=0$ for any $0 \leq \alpha \leq n_j$ and any $a
\in \ZZ$ unless $\alpha=0$ and $a \geq 0$ or $\alpha=n_j$ and $a
\leq -n_j-1$, together with the fact that $\cE_{n_1 + \cdots +n_s
}= \cO_{X}$ we deduce that for any  $t \in \ZZ$ and any $E_i^k \in
\cE_k$, $0 \leq k \leq n_1 + \cdots+ n_s-1$, $0 < \alpha < n_1 +
\cdots+ n_s$
\[ H^{\alpha}(X, \cO_{X}(t, \cdots, t)\otimes E_i^k)=0.\]
Hence, $X$ has the strong CM property.

(4) Let $X=Gr(k,n)$ be the Grassmannian variety of $k$-dimensional
subspaces of the $n$-dimensional vector space and take
$\sigma=(\cE_0, \cdots, \cE_{k(n-k)})$ be the $k(n-k)$-block
collection of vector bundles on $X$ given in  Example
\ref{exempleblock} (1). Notice that $\cE_{k(n-k)}=\cO_X$. Hence,
since any $\Sigma^{\alpha}S \in \cE_r$, $0 \leq r \leq k(n-k)-1$,
is an ACM vector bundle, we get that $X=Gr(k,n)$ has the CM
property but not the strong CM property.

(5) Let $\pi: \widetilde{\PP}^2(3) \rightarrow \PP^2$  be the blow
up of $\PP^2$ at $3$ points and let $L_i=\pi^{-1}(p_i)$, $1\leq i
\leq 3$, be the exceptional divisors. Then,
\[ (\cO, \cO(H), \cO(2H-L_1-L_2-L_3), \cO(2H-L_2-L_3), \cO(2H-L_1-L_3), \cO(2H-L_1-L_2)) \]
 is a full exceptional collection of coherent sheaves on
 $\widetilde{\PP}^2(3)$.   By \cite{KN}; Proposition 4.2 (3), the collection
  $(\cE_0, \cE_1, \cE_2)$ with $\cE_0=(\cO)$,
 $\cE_1= (\cO(H), \cO(2H-L_1-L_2-L_3))$ and  $\cE_2=(\cO(2H-L_2-L_3), \cO(2H-L_1-L_3),
 \cO(2H-L_1-L_2))$,  is a $3$-block collection
 of line bundles on $\widetilde{\PP}^2(3)$. Hence,
 $\widetilde{\PP}^2(3)$ has the weak CM property.
\end{example}

We are led to pose the following problem/question:

\begin{problem} To characterize smooth projective varieties with
the (weak, strong) CM property.
\end{problem}

By \cite{BP}; Theorem 3.4, any smooth projective variety with an
excellent collection is Fano. All examples described above about
smooth projective varieties with the (weak, strong) CM property
are Fano. So,  we wonder

\begin{question} Let $X$ be a smooth projective variety and assume
that X has the (weak, strong) CM property. Is X Fano?
\end{question}

\vspace{4mm} Beilinson's Theorem was stated in 1978 \cite{Be} and
since then it has became a major tool in classifying vector
bundles over projective spaces. Beilinson's spectral sequence was
generalized by Kapranov to hyperquadrics and Grassmannians
(\cite{Kap} and \cite{Kap2}) and by the authors to any smooth
projective variety with an excellent collection \cite{CMR3}. We
are now ready to generalize Beilinson's Theorem to any smooth
projective variety which has the weak CM property and to state the
main result of this section.

 \vspace{4mm}

\begin{theorem} \label{mbe} {\bf (Beilinson type spectral sequence)}
Let $X$ be a smooth projective  variety of dimension $n$
 with an $n$-block collection $\sigma =(\cE_0, \cE_1,\cdots ,\cE _{n})$,
  $\cE_i=(E_1^i, \ldots,E_{\alpha_i}^i )$ of coherent sheaves on $X$ which generates $\cD$.
   Then for any coherent sheaf $F$ on $X$ there are two spectral sequences with $E_1$-term
 \begin{equation}\label{sucespectral0}
   _{I} E_1^{pq}=
\begin{cases}
   \bigoplus_{i=1}^{\alpha_{p+n}}\Ext^q(R_{\cE_{n} \cdots \cE_{p+n+1}}E _{i}^{p+n},{F})
   \otimes E_{i}^{p+n}   & \mbox{if} \quad -n \leq p \leq -1 \\
\bigoplus_{i=1}^{\alpha_{n}}\Ext^q(E _{i}^{n},{F})
   \otimes E_{i}^{n}   & \mbox{if} \quad p=0 \end{cases}
   \end{equation}
 \begin{equation}\label{sucespectral1}
   _{II} E_1^{pq}= \begin{cases} \bigoplus_{i=1}^{\alpha_{p+n}}\Ext^q( (E_{i}^{p+n})^*,{F})
   \otimes (R_{\cE_{n} \cdots \cE_{p+n+1}}E _{i}^{p+n})^* &
   \mbox{if} \quad -n \leq p \leq -1 \\
\bigoplus_{i=1}^{\alpha_{n}}\Ext^q( {E_{i}^{n}}^*,{F})
   \otimes {E _{i}^{n}}^* &
   \mbox{if} \quad  p=0 \end{cases}  \end{equation}
    situated in the square $-n \le p  \le 0$, $0 \le q \le n $
   which converge to $$ _{I} E^{i}_{\infty }=_{II} E^{i}_{\infty
   }=\begin{cases} {F} \mbox{ for } i=0 \\ 0 \mbox{ for } i\ne 0.\end{cases}$$
\end{theorem}
\begin{proof}
We will only prove the existence of the first spectral sequence.
The other can be done similarly. For any $\gamma$, $0 \leq \gamma
\leq n$, we write $^{i}V^{\bullet}_\gamma$ for the graded vector
spaces
\[^{i}V^{\bullet}_\gamma = \Hom^{\bullet}_{\cD}
(R_{\cE_{n}\cdots \cE_{\gamma+1}}E_{i}^{\gamma}, F)=
\Hom^{\bullet}_{\cD}(E_{i}^{ \gamma}, L_{\cE_{\gamma+1}\cdots
\cE_n}F)\] where the second equality follows from standard
properties of mutations (\cite{GK}; Pag. 12-14).

By Remark \ref{triangles}, the triangles defining the consequent
right mutations of $F$ and the consequent left mutations of $F[n]$
through $(\cE_0, \cdots, \cE_n)$ can be written as
\[ (\bigoplus_{i=1}^{\alpha_ {\gamma}} {}^{i}V_{\gamma}^{\bullet} \otimes
E_i^{\gamma})[-1] \mapright{k_{\gamma}} R_{\cE_{\gamma}\cdots
\cE_{0}} F[-1] \mapright{i_{\gamma}} R_{\cE_{\gamma-1}\cdots
\cE_{0}} F \mapright{j_{\gamma}} \bigoplus_{i=1}^{\alpha_
{\gamma}} {}^{i}V_{\gamma}^{\bullet} \otimes E_i^{\gamma}\]
\[ \bigoplus_{i=1}^{\alpha_ {\gamma}} {}^{i}V_{\gamma}^{\bullet} \otimes
E_i^{\gamma} \mapright{j^{\gamma+1}} L_{\cE_{\gamma+1}\cdots
\cE_{n}} F[n] \mapright{i^{\gamma+1}} L_{\cE_{\gamma}\cdots
\cE_{n}} F[n+1] \mapright{k^{\gamma+1}}( \bigoplus_{i=1}^{\alpha_
{\gamma}} {}^{i}V_{\gamma}^{\bullet} \otimes E_i^{\gamma})[1].
\]
We arrange them into the following big diagram:

\vskip 6mm

 \xymatrix{
0=R_{\cE_n\cdots \cE_0}F \ar[dd]_{i_n} & & F[n] \ar[dd]^{i^0}
\\
 & \ar[ul]_{k_n} \bigoplus_{i=1}^{\alpha_n} {}^iV_n^{\bullet} \otimes E_i ^n
   \ar[ur]^{j^0} & \\
  R_{\cE_{n-1}\cdots \cE_0}F \ar[ur]^{j_n}\ar[dd]_{i_{n-1}} & &
  \ar[ul]_{k^0} \ar[dd]^{i^1} L_{\cE_n}F[n] \\
  & \ar[ul]_{k_{n-1}} \ar[uu]^{d_{n-1}}
   \bigoplus_{i=1}^{\alpha_{n-1}} {}^iV_{n-1}^{\bullet} \otimes E_i ^{n-1}  \ar[ur]^{j^1} & \\
R_{\cE_{n-2} \cdots \cE_0}F \ar[ur]^{j_{n-1}} \ar@{.}[ddd] & &
  \ar[ul]_{k^1} \ar@{.}[ddd] L_{\cE_{n-1} \cE_n}F[n] \\
  & \ar[ul]_{k_{n-2}} \ar[uu]^{d_{n-2}}
   \bigoplus_{i=1}^{\alpha_{n-2}} {}^iV_{n-2}^{\bullet} \otimes E_i ^{n-2}  \ar[ur]^{j^2} & \\
   && \\
   R_{\cE_{1} \cE_0}F  \ar[dd]_{i_1} & &
   \ar[dd]^{i^{n-1}} L_{\cE_{2} \cdots \cE_n}F[n] \\
  & \ar[ul]_{k_{1}} \ar@{.}[uuu]
   \bigoplus_{i=1}^{\alpha_{1}} {}^iV_{1}^{\bullet} \otimes E_i ^{1}  \ar[ur]^{j^{n-1}} & \\
R_{\cE_0}F \ar[ur]^{j_{1}} \ar[dd]_{i_0} & &
  \ar[ul]_{k^{n-1}} \ar[dd]^{i^n} L_{\cE_{1} \cdots \cE_n}F[n] \\
  & \ar[ul]_{k_{0}} \ar[uu]^{d_{0}}
   \bigoplus_{i=1}^{\alpha_{0}} {}^iV_{0}^{\bullet} \otimes E_i ^{0}  \ar[ur]^{j^n} & \\
 F \ar[ur]^{j_0} & & \ar[ul]_{k^{n}} L_{\cE_{0} \cdots \cE_n}F[n] =0
 \\}

At this diagram, all oriented triangles along left and right
vertical borders are distinguished, the morphisms  $i_{\bullet}$
and $i^{\bullet}$ have degree one, and all triangles and rhombuses
in the central column are commutative. So, there is the following
complex, functorial on $F$,

\[L^{\bullet}: 0 \rightarrow
\bigoplus_{i=1}^{\alpha_0}{}^{i}V_0^{\bullet} \otimes E_{i}^0
 \rightarrow \bigoplus_{i=1}^{\alpha_1}{}^{i}V_1^{\bullet} \otimes E_{i}^1 \rightarrow \cdots
 \rightarrow \bigoplus_{i=1}^{\alpha_{n-1}}{}^{i}V_{n-1}^{\bullet} \otimes
 E_{i}^{n-1}
  \rightarrow \bigoplus_{i=1}^{\alpha_n}{}^{i}V_{n}^{\bullet} \otimes
  E_{i}^{n}
 \rightarrow 0   \]
 and by the above Postnikov-system we have that $F$ is a right convolution of this complex.
  Then, for an
 arbitrary linear covariant cohomological functor
 $\Phi^{\bullet}$, there exists an spectral sequence with
 $E_1$-term
 \[ _{I}E_1^{pq}=\Phi^q(L^p)\]
 situated in the square $0 \leq p,q \leq n$ and converging to
  $\Phi^{p+q}(F)$ (see \cite{Kap2}; 1.5). Since
 $\Phi^{\bullet}$ is a linear functor, we have
\begin{equation} \label{espectralgral}
 \Phi^q(L^p)= \bigoplus_{i=1}^{\alpha_p}\Phi^q(^{i}V_p^{\bullet} \otimes E_{i}^{p})=
\bigoplus_{i=1}^{\alpha_p}\bigoplus_{l} {}^{i}V_p^l \otimes
\Phi^{q-l}(E_{i}^{p})= \bigoplus_{i=1}^{\alpha_p}
\bigoplus_{\alpha+\beta=q}  {}^{i}V_p^{\alpha} \otimes
\Phi^{\beta}(E_{i}^{p}). \end{equation} In particular, if we
consider the covariant linear cohomology functor which takes a
complex to its cohomology sheaf and acts identically on pure
sheaves, i.e.
\[ \Phi^{\beta}(F)= \begin{cases} F \mbox{ for } \beta=0 \\ 0 \mbox{ for } \beta \neq 0
 \end{cases}\]
 on any pure sheaf $F$, in the square $0 \leq p,q \leq n$, we get
 \[ _{I} E_1^{pq}=\bigoplus_{i=1}^{\alpha_p}
 {}^{i}V_p^{q} \otimes E_{i}^{p}= \bigoplus_{i=1}^{\alpha_p}
 \Ext^q(R_{\cE_{n} \cdots \cE_{p+1}}E _{i}^p,{F})
   \otimes E_{i}^p \]
 which converges to
 $$ _{I} E^{i}_{\infty }=
\begin{cases} {F} \mbox{ for } i=0 \\ 0 \mbox{ for } i\ne
0.\end{cases}$$ Finally, if we call $p'=p-n$, we get the spectral
sequence $$   _{I}
E_1^{p'q}=\bigoplus_{i=1}^{\alpha_{p'+n}}\Ext^q(R_{\cE_{n} \cdots
\cE_{p'+n+1}}E _{i}^{p'+n},{F})
   \otimes E_{i}^{p'+n}$$
  situated in the square $-n \le p'  \le 0$, $0 \le q \le n $
   which converges to $$ _{I} E^{i}_{\infty }=\begin{cases} {F} \mbox{ for } i=0 \\ 0 \mbox{ for } i\ne 0.\end{cases}$$
\end{proof}


\section{Splitting vector bundles and cohomological
characterization of vector bundles}

A well known result of  Horrocks states that a vector bundle on
$\PP^n$ has no intermediate cohomology if and only if it splits
into a direct sum of line bundles. The first  goal of this section
is to generalize Horrocks criterion to vector bundles on
multiprojective spaces  and to  any smooth projective variety with
the strong CM property. As a main tool we will use the
Beilinson-type spectral sequences stated in the previous section.

\begin{proposition}\label{main1}
Let $X$ be a smooth projective  variety of dimension $n$
 with the CM property given by the $n$-block collection $\sigma =(\cE_0, \cE_1,\cdots ,\cE _{n})$,
  $\cE_i=(E_1^i, \ldots,E_{\alpha_i}^i )$ of coherent sheaves on $X$. Let $F$
   be a coherent sheaf on $X$ such that for any $-n  \le p \le
   -1
   $ and $1 \le i \leq \alpha_p$
   \[ H^{-p-1}(X, F \otimes E_i^{p+n})=0 . \]
Then $F$ contains
$\bigoplus_{i=1}^{\alpha_n}({E_i^n}^*)^{h^0(F\otimes E_i^n )}$ as
a direct summand.
\end{proposition}
\begin{proof}
By Theorem \ref{mbe}, there is a spectral sequence with $E_1$-term

$$    _{II} E_1^{pq}= \begin{cases}
\bigoplus_{i=1}^{\alpha_{p+n}}\Ext^q( (E_{i}^{p+n})^*,{F})
   \otimes (R_{\cE_{n} \cdots \cE_{p+n+1}}E _{i}^{p+n})^* &
   \mbox{if} \quad -n \leq p \leq -1 \\
\bigoplus_{i=1}^{\alpha_{n}}\Ext^q( {E_{i}^{n}}^*,{F})
   \otimes {E _{i}^{n}}^* &
   \mbox{if} \quad  p=0 \end{cases} $$
    situated in the square $-n \le p  \le 0$, $0 \le q \le n $
   which converges to $$ _{II} E^{i}_{\infty }
   =\begin{cases} {F} \mbox{ for } i=0 \\ 0 \mbox{ for } i\ne 0.\end{cases}$$
   By assumption, $ _{II} E_1^{p,-p-1}=0$, i.e., the $E_1$-term
   looks like

$$ 
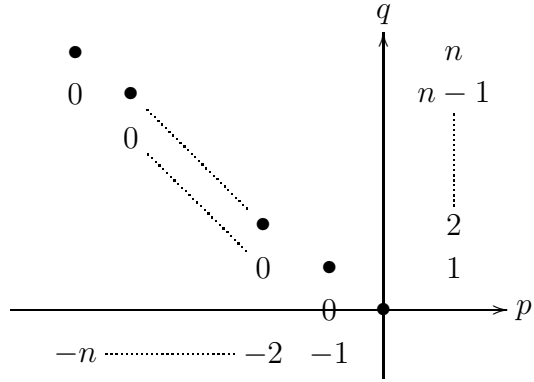
\xymatrix@R=0.05cm@C=0.1cm{ & & & & & & & & & q & & \\  & &
\bullet &&&&& & && n & \\ & & 0 &\bullet & & & & & && n-1 & \\ & &
& 0 & & & & & && &
\\ && & & & & & & && & \\ && &
&& &  & & &&  & \\ & & & & & & & \bullet \ar@{.}[uuuullll] & && 2
\ar@{.}[uuuu]
 & \\ & & & & & & & 0 \ar@{.}[uuuullll] & \bullet& & 1 &
\\ \ar[rrrrrrrrrrr] & & & & &  &  &   &  0  &   \bullet
 &  & p \\ & & -n &  &  &  & & -2 \ar@{.}[lllll] &
 -1&   & &   \\ &&&&&&&&&   \ar[uuuuuuuuuu]  & &} $$

   So, the limit $ _{II} E^{i}_{\infty }$, i.e., $F$, contains
   $_{II} E_1^{00}=\bigoplus_{i=1}^{\alpha_n}({E_i^n}^*)^{h^0(F\otimes
E_i^n )}$ as a direct summand.
\end{proof}

As an immediate consequence of Proposition \ref{main1} we will
first re-prove Horrocks criterion.

\begin{corollary} \label{coro1}
Let $E$ be a vector bundle on $\PP^n$. The following conditions
are equivalent:
\begin{itemize}
\item[(i)] $E$ splits into a sum of line bundles.
 \item[(ii)] $E$
has no intermediate cohomology; i.e. $H^{i}(\PP^n,E(t))=0$ for
$1\le i \le n-1$ and for all $t\in \ZZ$
\end{itemize}
\end{corollary}
\begin{proof} $(i) \Rightarrow (ii)$. It follows from Bott's
formula. $(ii) \Rightarrow (i)$. We may suppose that $E$ is
indecomposable. So that it suffices to prove that $E$ is a line
bundle. To this end, we choose an integer $m$ such that
$H^0(\PP^n,E(m-1))=0$ and $H^0(\PP^n,E(m)) \neq 0$ and we apply
Proposition \ref{main1} to $X =\PP^n$, $\sigma=(\cO_{\PP^n}(-n),
\cdots, \cO_{\PP^n}(-1), \cO_{\PP^n})$ and $F=E(m)$. We conclude
that $\cO^{h^0E(m)}$ is a direct summand of $F$ and since $F$ is
indecomposable we get that $F=\cO_{\PP^n}$ and we are done.
\end{proof}

In \cite{O2} Ottaviani pointed out that Horrocks criterion fails
on a non singular quadric hypersurface $Q_n\subset \PP^{n+1}$; the
Spinor bundles $S$ on $Q_n$ have no intermediate cohomology and
they do not decompose into a direct sum of line bundles.
Nevertheless, we have the following cohomological characterization
of vector bundles on $Q_n$ which split into a direct sum of line
bundles; and of vector bundles on a Grassmannian $Gr(k,n)$ which
also split into a direct sum of line bundles.

\vspace{3mm}

On $Q_n$, we shall use the unified notation $\Sigma_{\ast}$
meaning that for even $n$ both Spinor bundles $\Sigma_1$ and
$\Sigma_2$ are considered, and for odd $n$, the Spinor bundle
$\Sigma$ (see Example  \ref{exempleblock} (2) for more details).

\begin{corollary} \label{coro2} Let $E$ be a vector bundle on $Q_n\subset
\PP^{n+1}$. The following conditions are equivalent:
\begin{itemize}
\item[(i)] $E$ splits into a sum of line bundles.
 \item[(ii)]
$H^{i}(Q_n,E(t))=0$ for $1\le i \le n-1$ and  $t\in \ZZ$; and
$H^{n-1}(Q_n,E\otimes \Sigma_{\ast}(t-n))=0$.
\end{itemize}
\end{corollary}
\begin{proof}  $(i) \Rightarrow (ii)$. It is a well known statement. $(ii) \Rightarrow
(i)$. We may suppose that $E$ is indecomposable. So that it
suffices to prove that $E$ is a line bundle. To this end, we
choose an integer $m$ such that $H^0(Q_n,E(m-1))=0$ and
$H^0(Q_n,E(m)) \neq 0$ and we apply Proposition \ref{main1} to $X
=Q_n$, $\sigma=(\cE_0, \cdots, \cE_{n})$ defined  in Example
\ref{exempleblock} (2) and $F=E(m)$ (see also Example
\ref{CMexample}). Hence, we obtain that $\cO_{Q_n}^{h^0E(m)}$ is a
direct summand of $F$ and since $F$ is indecomposable we conclude
that $F=\cO_{Q_n}$.
\end{proof}

\vspace{3mm}

Keeping the notations introduced in Example \ref{exempleblock}
(1), we have: \vspace{3mm}

\begin{corollary} \label{coro3} Let $E$ be a vector bundle on
$Gr(k,n)$ and set $$\cE_r=\{\Sigma^{\alpha}S  | k(n-k)-r
=|\alpha|\}.$$ The following conditions are equivalent:
\begin{itemize}
\item[(i)] $E$ splits into a sum of line bundles.
 \item[(ii)] $H^{i}(Gr(k,n),E(t) \otimes \Sigma^{\alpha}S)=0$ for
$1\le i \le k(n-k)-1$, $t \in \ZZ$ and $\Sigma^{\alpha} S \in
\cE_{k(n-k)-i-1}$.
\end{itemize}
\end{corollary}
\begin{proof} $(i) \Rightarrow (ii)$. It is a well known statement. $(ii) \Rightarrow (i)$.
We may suppose that $E$ is indecomposable. So that it suffices to
prove that $E$ is a line bundle. To this end, we choose an integer
$m$ such that $H^0(Gr(k,n),E(m-1))=0$ and $H^0(Gr(k,n),E(m)) \neq
0$. We consider Proposition \ref{main1} applied  to $X =Gr(k,n)$,
$\sigma=(\cE_0, \cdots,  \cE_{k(n-k)})$ given in Example
\ref{exempleblock} (1) and $F=E(m)$ (see also Example
\ref{CMexample}) and we get that $\cO_{Gr(k,n)}^{h^0E(m)}$ is a
direct summand of $F$. Since $F$ is indecomposable we derive that
$F=\cO_{Gr(k,n)}$ and we are done.
\end{proof}

\begin{remark} Applying again Proposition \ref{main1} and arguing
as in Corollaries \ref{coro2} and \ref{coro3}, we can deduce the
splitting criteria for vector bundles on the Fano 3-folds  $V_5$
and $V_{22}$ given by Faenzi in \cite{Fa} and \cite{Fa1}.
\end{remark}

\vspace{2mm}
\begin{theorem} \label{criterion}
Let $X$ be a smooth projective variety of dimension $n$ with the
strong CM property given by the $n$-block collection
$\sigma=(\cE_0, \cdots, \cE_n)$, $\cE_i=(E_1^i,
\ldots,E_{\alpha_i}^i )$, of line bundles on $X$. Let $E$ be a
vector bundle on $X$ such that $E \otimes E_j^i$ is an ACM bundle
for any  $E_j^i \in \cE_i$, $0 \leq i \leq n-1$. Then, $E$ splits
into a direct sum of line bundles.
\end{theorem}
\begin{proof} We may suppose that $E$ is
indecomposable. So that it suffices to prove that $E$ is a line
bundle. By assumption, for any  $E_j^i \in \cE_i$, $0 \leq i \leq
n-1$, any $0 < p < n$ and any $t \in \ZZ$,
\[ H^p(X,E \otimes E_j^i \otimes \cO_X(t) )=0. \]
  We choose an integer $m$ such that
  $$\oplus _{j=1}^{\alpha _n}H^0(X,E\otimes \cO_X(m-1)\otimes
  E_j^{n})=0\text{ and }
  \oplus _{j=1}^{\alpha _n}H^0(X,E\otimes \cO_X(m)\otimes E_j^{n})\ne 0.$$
We apply Proposition \ref{main1} to $X$, $\sigma =(\cE_0,
\cE_1,\cdots ,\cE _{n})$ and $F=E(m)$. We conclude that $F$
contains $\bigoplus_{i=1}^{\alpha_n}({E_i^n}^*)^{h^0(F\otimes
E_i^n )}$ as a direct summand and since $F$ is indecomposable we
get that $F={E_i^n}^*$ for some $1\le i \le \alpha _n$ which
proves what we want.
\end{proof}

As a consequence we get:

\vspace{3mm}

\begin{theorem} \label{mainmulti} Let  $X=\PP^{n_1}\times \cdots \times \PP^{n_r}$
be a multiprojective space and let $E$ be a vector bundle on $X$
such that  $E\otimes \cO_X(t_1,\cdots ,t_r)$ is an ACM bundle for
any $-n_i \leq t_i \leq 0$, $1 \leq i \leq r$. Then, $E$ splits
into a direct sum of line bundles.
\end{theorem}
\begin{proof} Let $\sigma=(\cE_0, \cdots, \cE_{n_1 +
\cdots + n_r })$ be the $(n_1+\cdots+n_r)$-block collection of
line bundles on $X$ given in Example \ref{exempleblock} (3) (see
also Example \ref{CMexample}). Then, we apply Theorem
\ref{criterion}.
\end{proof}

The converse of Theorem \ref{criterion} turns to be true for
vector bundles on projective spaces (Horrock's criterion) but, in
general, it is not true.  For instance, as a consequence of the
K\"{u}nneth formula,  on any multiprojective space
$\PP^{n_1}\times \cdots \times \PP^{n_r}$ there are many line
bundles $L$ such that $L \otimes \cO(t_1, \cdots, t_r)$ is not an
ACM bundle (take for example $L=\cO_{\PP^2 \times \PP^3}(-3,4)$).

\vspace{2mm}

As another application of Beilinson-type spectral sequence we will
derive a cohomological characterization of huge families of vector
bundles. The first attempt in this direction is due to Horrocks
who in \cite{Ho2} gave a cohomological characterization of the
sheaf of $p$-differential forms, $\Omega ^p _{\PP^n}$. Similarly,
in \cite{AO}, Ancona and Ottaviani obtained a cohomological
characterization of the vector bundles $\psi _{i}$ on $Q_n$
introduced by Kapranov in \cite{Kap2}. These two results are a
particular case of this following much more general statement.

\begin{proposition}\label{main2}
Let $X$ be a smooth projective  variety of dimension $n$
 with an $n$-block collection $\sigma =(\cE_0, \cE_1,\cdots ,\cE _{n})$,
  $\cE_i=(E_1^i, \ldots,E_{\alpha_i}^i )$ of coherent sheaves on $X$ which generates $\cD$
   and let $F$ be a coherent sheaf on $X$. Assume there exists $j$,
   $0 < j < n$ such that for any $-n  \le p \le -j-1
   $ and $1 \le i \leq \alpha_p$
   \[ H^{-p-1}(X, F \otimes E_i^{p+n})=0  \]
   and for any   $-j+1  \le p \le 0
   $ and $1 \le i \leq \alpha_p$
   \[ H^{-p+1}(X, F \otimes E_i^{p+n})=0 . \]
Then $F$ contains $\bigoplus_{i=1}^{\alpha_{n-j}} ((R_{\cE_{n}
\cdots \cE_{n+1-j}}E _{i}^{n-j})^*)^{h^j(F\otimes E_i^{n-j})}$ as
a direct summand.
\end{proposition}
\begin{proof}
By Theorem \ref{mbe}, there is a spectral sequence with $E_1$-term
$$    _{II} E_1^{pq}= \begin{cases}
\bigoplus_{i=1}^{\alpha_{p+n}}\Ext^q( (E_{i}^{p+n})^*,{F})
   \otimes (R_{\cE_{n} \cdots \cE_{p+n+1}}E _{i}^{p+n})^* &
   \mbox{if} \quad -n \leq p \leq -1 \\
\bigoplus_{i=1}^{\alpha_{n}}\Ext^q( {E_{i}^{n}}^*,{F})
   \otimes {E _{i}^{n}}^* &
   \mbox{if} \quad  p=0 \end{cases} $$
    situated in the square $-n \le p  \le 0$, $0 \le q \le n $
   which converges to $$ _{II} E^{i}_{\infty }
   =\begin{cases} {F} \mbox{ for } i=0 \\ 0 \mbox{ for } i\ne 0.\end{cases}$$
   By assumption, there exists an integer $j$, $0<j<n$ such that
    $_{II} E_1^{p,-p-1}=0$ for  any $-n  \le p \le -j-1
   $ and  $_{II} E_1^{p,-p+1}=0$ for  any $-j+1 \leq p \le 0$. Therefore, we have the
   following $E_1$-diagram

$$ 
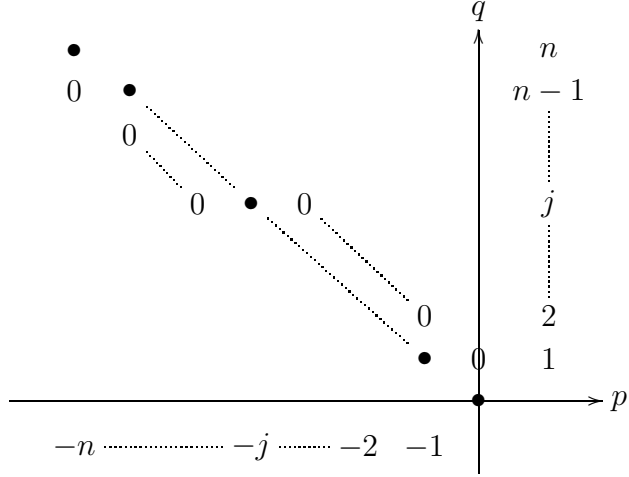
\xymatrix@R=0.05cm@C=0.1cm{ & & & & & & & & & & q & & \\  & &
\bullet &&&&& & & && n & \\ & & 0 &\bullet & & & & & & && n-1 & \\
& & & 0 & & & && & && &
\\ && & & & & && & && & \\
&& & & &  0 \ar@{.}[uull]& \bullet \ar@{.}[uuulll] & 0
\ar@{.}[ddddrr]& && & j \ar@{.}[uuu]&
\\ && & && & &&  & && &
\\ && & && & &&  & && &
\\ && & && & &&  & && &
\\ & & & & & & & &   &  0&& 2 \ar@{.}[uuuu]
 & \\ & & & & & & &&  &  \ar@{.}[uuuuulll] \bullet& 0& 1 &
\\ \ar[rrrrrrrrrrrr] &&
 & &
&
 &  & &
  &    &   \bullet
 &  & p \\ & & -n \ar@{.}[rrrr]&  &  &  &   -j & & -2 \ar@{.}[ll] &
 -1&   & &   \\ & &&&&&&&&&   \ar[uuuuuuuuuuuuu]  & & }$$

   So, the vector bundle  $F$ contains
   $_{II} E_1^{jj}=\bigoplus_{i=1}^{\alpha_{n-j}} ((R_{\cE_{n}
\cdots \cE_{n+1-j}}E _{i}^{n-j})^*)^{h^j(F\otimes E_i^{n-j})}$ as
a direct summand.
\end{proof}

\vspace{3mm}

Our next goal is to extend Horrock's characterization of
$p$-differentials  over $\PP^n$ to multiprojective spaces
$\PP^{n_1} \times \cdots \times \PP^{n_s}$. To this end, we will
first determine the left dual $(n_1+ \cdots + n_s)$-block
collection of the $(n_1+ \cdots + n_s)$-block collection
$\sigma=(\cE_0, \cdots, \cE_{(n_1+ \cdots + n_s)})$ described in
Example \ref{exempleblock}.

\vspace{3mm}

\begin{notation}
Let $X_1$ and $X_2$ be two smooth projective varieties and let
$$p_i: X_1 \times X_2 \rightarrow X_i, \quad  i=1,2,$$ be the
natural projections. We denote by $B_1 \boxtimes B_2$ the exterior
tensor product of $B_i$ in $\cO_{X_i}$-mod, $i=1,2$, i.e. $B_1
\boxtimes B_2=p_1^*B_1 \otimes p_2^*B_2$ in $\cO_{X_1 \times
X_2}$-mod.
\end{notation}

\vspace{3mm}

\begin{proposition}
\label{erredemulti} Let $X=\PP^{n_1} \times \cdots \times
\PP^{n_s}$ be a multiprojective space of dimension
$d=n_1+\cdots+n_s$.
 For any $0 \leq j \leq d$, denote by $\cE_j$
 the collection of all line bundles on $X$
 \[ \cO_X(a_1^j, a_2^j, \cdots, a_s^j)\]
 with $-n_i \leq a_i^j \leq 0$ and
 $\sum_{i=1}^{s}a_i^j=j-d$. Then,
  for any $\cO_{X}(t_1,\cdots,t_s) \in \cE_{d-k}$ and  any $0 \leq k \leq d$,
\[R^{(k)}\cO_{X}(t_1,\cdots,t_s)=
 R_{\cE_d \cdots \cE_{d-k+1}}\cO_{X}(t_1,\cdots,t_s)=
  \bigwedge^{-t_1} T_{\PP^{n_1}}(t_1)
\boxtimes \cdots \boxtimes \bigwedge^{-t_s} T_{\PP^{n_s}}(t_s).\]
\end{proposition}
\begin{proof} According to Proposition \ref{dual}
(\ref{caracdelsr}), we only need to see that $\bigwedge^{-t_1}
T_{\PP^{n_1}}(t_1) \boxtimes \cdots \boxtimes \bigwedge^{-t_s}
T_{\PP^{n_s}}(t_s)$ verifies the orthogonality conditions
(\ref{orto1}) and (\ref{orto2}). For any $i$, $0 \leq i \leq d $,
let $\cO_{X}(a_1^i, \cdots , a_s^i) \in \cE_ i$. By the
K\"{u}nneth formula,
\[ H^{\alpha}(X,\bigwedge^{-t_1}
\Omega_{\PP^{n_1}}(-t_1) \boxtimes \cdots \boxtimes
\bigwedge^{-t_s} \Omega_{\PP^{n_s}}(-t_s) \otimes \cO_{X}(a_1^i,
\cdots , a_s^i))  \] \[=
\bigoplus_{\alpha_1+\cdots+\alpha_s=\alpha}
H^{\alpha_1}(\PP^{n_1},\bigwedge^{-t_1}\Omega_{\PP^{n_1}}(a_1^i-t_1))
\otimes \cdots \otimes H^{\alpha_s}(\PP^{n_s}, \bigwedge^{-t_s}
\Omega_{\PP^{n_s}}(a_s^i-t_s)). \] Using Bott's formula, it is
zero unless $\alpha=k$, $i=d-k$ and $\cO_{X}(a_1^i, \cdots ,
a_s^i)=\cO_{X}(t_1,\cdots,t_s)$, which proves what we want.
\end{proof}

\vspace{3mm}

The following result gives us a precise cohomological
characterization of sheaves of $p$-differential forms on
multiprojective spaces.

\vspace{3mm}

\begin{theorem}
\label{mainmulti2}
 Let $X=\PP^{n_1} \times \cdots \times
\PP^{n_s}$ be a multiprojective space of dimension
$d=n_1+\cdots+n_s$.
 For any $0 \leq i \leq d$, denote by
  $\cE_i=(E_1^i, \cdots, E_{\alpha_i}^i)$
 the collection of all line bundles on $X$
 \[ \cO_X(a_1^i, a_2^i, \cdots, a_s^i)\]
 with $-n_k \leq a_k^i \leq 0$ and
 $\sum_{k=1}^{s}a_k^i=i-d$. Assume there exists a rank ${d \choose
 j}$vector bundle $F$ on $X$ with
   $0 < j < d$, such that for any $-d  \le p \le -j-1
   $ and $1 \le i \leq \alpha_p$
   \[ H^{-p-1}(X, F \otimes E_i^{p+d})=0,  \]
   for any   $-j+1  \le p \le 0   $ and $1 \le i \leq \alpha_p$
   \[ H^{-p+1}(X, F \otimes E_i^{p+d})=0  \]
   and $H^j(F\otimes E_i^{d-j})= \CC$ for any $1 \leq i \leq \alpha_{d-j}$.
    Then $F$ is isomorphic to the bundle of $(d-j)$-differential forms,
    i.e.
   \[ F \cong \bigwedge^{d-j} (\Omega_{\PP^{n_1} \times  \cdots \times \PP^{n_s}}(1,\cdots,1))
    \cong \bigoplus_{t_1+\cdots+t_s=j-d}\bigwedge^{-t_1}
\Omega_{\PP^{n_1}}(-t_1) \boxtimes \cdots \boxtimes
\bigwedge^{-t_s} \Omega_{\PP^{n_s}}(-t_s) \] being
$E_i^{d-j}=\cO_X(t_1, \cdots,t_s)$.
\end{theorem}
\begin{proof} It follows from Proposition \ref{main2} and
Proposition \ref{erredemulti}.
\end{proof}

\vspace{3mm}

We will end this section extending Horrocks characterization of
sheaves of $p$-differential forms in $\PP^n$ to Grassmannians.
Notice that under the isomorphism $Gr(1,n+1) \cong \PP^{n}$, the
universal quotient bundle $\cQ$ on $Gr(1,n+1)$ corresponds to
$\Omega_{\PP^n}(1)$. So, it is natural to get, as a generalization
of Horrocks characterization of the bundles
$\Omega^p_{\PP^n}(p)=\bigwedge^p (\Omega_{\PP^n}(1))$, a
cohomological characterization of the bundles $\Sigma^{\beta}\cQ$
being $\cQ$ the universal quotient bundle on $Gr(k,n)$. More
precisely, keeping the notations introduced in Example
\ref{prihirse} (1) and in Example \ref{exempleblock} (1) we have:

\vspace{3mm} According to Example \ref{prihirse}, for any
$\beta=(\beta_1, \cdots, \beta_{n-k})$ with $k \geq \beta_1 \geq
\beta_2 \geq \cdots \geq \beta_{n-k} \geq 1$, denote by
$r_{\beta}$ the rank of $\Sigma^{\beta}\cQ$ and consider
$r_j=\sum_{|\beta|=j} r_{\beta}$.

\vspace{3mm}

\begin{corollary}
Let $F$ be a vector bundle on $Gr(k,n)$, set $d=k(n-k)$ and
$$\cE_r=\{\Sigma^{\alpha}\cS  | k(n-k)-r =|\alpha|\}.$$ Assume
there exists $j$,   $0 < j < d$ such that for any $-d  \le p \le
-j-1
   $, $1 \le i \leq \alpha_p$ and any $\Sigma^{\alpha}\cS \in \cE_{d+p}$
   \[ H^{-p-1}(Gr(k,n), F \otimes \Sigma^{\alpha}\cS)=0  \]
   and for any   $-j+1  \le p \le 0
   $, $1 \le i \leq \alpha_p$ and any $\Sigma^{\alpha}\cS \in \cE_{d+p}$
   \[ H^{-p+1}(Gr(k,n), F \otimes \Sigma^{\alpha}\cS)=0 . \]
If  $\rank F= r_j$ then, $F$ is isomorphic to
$\bigoplus_{|\beta|=j} \Sigma^{\beta} \cQ^*$.
\end{corollary}
\begin{proof} It is well known that the following orthogonality
relation between the bundles $\Sigma^{\alpha}\cS$ and
$\Sigma^{\beta}\cQ^*$ holds:
\[ H^q(Gr(k,n), \Sigma^{\alpha}\cS \otimes \Sigma^{\beta}\cQ^*)=
\begin{cases} \CC & \mbox{if }
\alpha=\widetilde{\beta} \quad \mbox{and } q=|\alpha| \\ 0 &
\mbox{otherwise.}
\end{cases}\]
So, the bundles $\Sigma^{\beta}\cQ^*$ verify the orthogonality
conditions (\ref{orto1}) and (\ref{orto2}) and we apply
Proposition \ref{main2}.
\end{proof}

\section{Final comments}

In \cite{R} Rouquier introduced the notion of dimension for a
triangulated category and he determined bounds for the dimension
of the bounded derived category  $D^b({\cO}_X-mod)$ of coherent
sheaves over an algebraic variety $X$. In particular, among other
results, he proved that if the diagonal of an algebraic variety
$X$ has a resolution of length $r+1$ then $\dim D^b({\cO}_X-mod)
\leq r$ and for any $n$-dimensional smooth projective variety $X$
we have  $n\le \dim D^b({\cO}_X-mod) \leq 2n$ . He also posed the
following questions:

\vspace{3mm}

\begin{question}
\label{qq} Does the inequality
\[\dim D^b({\cO}_{X \times Y}-mod) \leq \dim D^b({\cO}_X-mod) +
\dim D^b({\cO}_Y-mod) \] hold for $X$, $Y$ separated schemes of
finite type over a perfect field?
\end{question}

\begin{question}
\label{qqq} Is there any example of $n$-dimensional smooth
projective variety $X$ with $n< \dim D^b({\cO}_X-mod)$?

\end{question}
\vspace{3mm}

Using the results we have obtained in this paper, we are able to
contribute to these questions and we will prove that the equality
in Question \ref{qq} holds for multiprojective spaces and we will
enlarge the family of $n$-dimensional smooth projective variety
$X$ such that  $n= \dim D^b({\cO}_X-mod) \leq 2n$. Indeed, we have

\vspace{3mm}
\begin{theorem}
\label{dim} Let $X$ be  a smooth projective variety with the weak
CM property. Then
\[ \dim D^b({\cO}_X-mod)= \dim X.\]
 \end{theorem}
 \begin{proof}
Denote by $n$ the dimension of $X$ and consider an $n$-block
collection $\sigma=(\cE_0, \cdots, \cE_n)$ with $\cE_i=(E_1^i,
\cdots, E_{\alpha_i}^i)$. Such $n$-block collection exists because
$X$ has the weak CM-property. By Theorem \ref{mbe}, we have the
following resolution of the diagonal
\[ 0 \rightarrow \bigoplus_{i=1}^{\alpha_0} (R_{\cE_{n} \cdots \cE_{1}}E
_{i}^{0})^*   \boxtimes E_{i}^{0} \rightarrow
\bigoplus_{i=1}^{\alpha_1} (R_{\cE_{n} \cdots \cE_{2}}E
_{i}^{1})^* \boxtimes E_{i}^{1} \rightarrow \cdots
\hspace{40mm}\]
\[ \hspace{30mm} \cdots \rightarrow \bigoplus_{i=1}^{\alpha_{n-1}} (R_{\cE_{n}}E
_{i}^{n-1})^*   \boxtimes E_{i}^{n-1} \rightarrow
\bigoplus_{i=1}^{\alpha_n} (E _{i}^{n})^* \boxtimes E_{i}^{n}
\rightarrow \cO_{\Delta} \rightarrow 0.  \] So, according to
\cite{R}; Proposition 5.5, $\dim D^b({\cO}_X-mod) \leq \dim X$. On
the other hand, by \cite{R}; Proposition 5.36, $\dim X \leq \dim
D^b({\cO}_X-mod)$     and we are done.
\end{proof}

In particular, we have:

\begin{proposition}
\label{dimmulti} Let $X=\PP^{n_1} \times \cdots \times \PP^{n_s}$
be a multiprojective space. Then
\[ \dim D^b({\cO}_{\PP^{n_1} \times \cdots \times
\PP^{n_s}}-mod)= \sum_{i=1}^{s} \dim D^b({\cO}_{\PP^{n_i}}-mod).\]
\end{proposition}
\begin{proof}
Since by Example \ref{CMproperty} (3), $X$ has the weak CM
property, the result follows from Theorem \ref{dim} and the fact
that, by \cite{R}; Example 5.6, $\dim
D^b({\cO}_{\PP^{n_i}}-mod)=n_i$ for any $1 \leq i \leq s$.
\end{proof}


\begin{thebibliography}{999}




 \bibitem{AC}E. Arrondo, L. Costa, {\em Vector
bundles on Fano 3-folds without intermediate cohomology}, Comm.
Alg., {\bf 28}  (2000), 3899-3911.

\bibitem{AG} E. Arrondo, B. Gra\~{n}a,  {\em Vector bundles on G(1,4)
without intermediate cohomology},
 Journal of Alg. {\bf 214}, (1999), 128-142.

\bibitem{AO} V. Ancona, G. Ottaviani {\em Some applications of
Beilinson's theorem to projective spaces and quadrics}, Forum
Math. {\bf 3} (1991), 157-176.

\bibitem{Be} A.A. Beilinson, {\em Coherent sheaves on $\PP^n$ and Problems of Linear Algebra},
Funkt. Anal. Appl. {\bf 12} (1979), 214-216.

 \bibitem{Bo} A.I. Bondal, {\em Representation of associative algebras
  and coherent sheaves}, Math. USSR Izvestiya {\bf 34} (1990),
  23-42.


\bibitem{BP} A.I. Bondal, A.E. Polishchuk, {\em Homological properties of associative
algebras: the method of helices}, Russian Acad. Sci. Izv. Math.
{\bf 42} (1994),  219-259.

\bibitem{BGS} R.O. Buchweitz, G.M. Greuel, F.O. Schreyer, {\em
Cohen-Macaulay modules on hypersurface singularities},
  Invent. Math.  {\bf 88}  (1987),  165-182.


\bibitem{CMR} L. Costa and R.M. Mir\'o-Roig, {\em Tilting sheaves
on toric varieties}, Math. Z., {\bf 248} (2004), 849-865.


\bibitem{CMR3} L. Costa and R.M. Mir\'o-Roig, {\em Tilting
bundles, helix theory and Castelnuovo-Mumford regularity},
Preprint 2004.



\bibitem{Fa} D. Faenzi, {\em Bundles over the Fano threefold $V_5$}, Communications in Algebra, to appear.



\bibitem{Fa1} D. Faenzi, {\em Bundles over  Fano threefolds of type $V_{22}$}, Annali di Matematica Pura e Applicata,
to appear.



\bibitem{GK} A.L. Gorodentsev, S.A. Kuleshov, {\em Helix Theory}, Mosc. Math. J.  4  (2004),  no. 2, 377-440, 535.


\bibitem{H} L. Hille {\em Consistent algebras and special tilting sequences},
  Math. Z.   {\bf 220} (1995), 189-205.


\bibitem{Ho} G. Horrocks, {\em Vector bundles on the punctured spectrum of
a local ring}, Proc. London Math. Soc. (3) {\bf 14} (1964),
689-713.

\bibitem{Ho2} G. Horrocks, {\em Construction of bundles on $\PP^n$}, Seminaire Douady-Verdier
ENS 77/78, Asterisque  {\bf 71-72} (1980), 197-203.


\bibitem{Kap} M. M. Kapranov, {\em On the derived category of coherent sheaves on Grassmann
manifolds}, Math. USSR Izvestiya, {\bf 24} (1985), 183-192.

\bibitem{Kap2} M. M. Kapranov, {\em On the derived category of coherent sheaves on
some homogeneous spaces}, Invent. Math., {\bf 92} (1988), 479-508.


\bibitem{KN} B. V. Karpov, D. Yu Nogin, {\em Three-block exceptional collections
over Del Pezzo surfaces}, Math. USSR Izvestiya, {\bf 62} (1998),
429-463.


\bibitem{O1} G. Ottaviani, {\em Crit\`{e}res de scindage pour les fibr\`{e}s
vectoriels sur les grassmanniens et les quadriques}, C.R. Acad.
Sci. {\bf 305} (1987), 257-260.

\bibitem{O2} G. Ottaviani, {\em Some extensions of Horrocks
criterion to vector bundles on grassmannians and quadrics}, Annali
di Matem. {\bf 155}, (1989), 317-341.


\bibitem{R} R. Rouquier, {\em Dimensions of triangulated categories},
 math.CT/0310134 (2003).


\end{thebibliography}
\end{document}